\documentclass{amsart}
\usepackage{color,graphicx,wrapfig}
\title[Uniform Regularity]{Uniform Regularity close to Cross Singularities
in an Unstable Free Boundary Problem} \dedicatory{Dedicated to Nina
Nikolaevna Uraltseva on the occasion of her 75th birthday}
\author[J. Andersson]{John Andersson}
\address{Department of Mathematics and Statistics, University of Jyv\"askyl\"a, Finland}
\email{johnande@jyu.fi}
\author[H. Shahgholian ]{Henrik Shahgholian}
\address{Department of Mathematics, Royal Institute of Technology,
100~44  Stockholm, Sweden}
\email{henriksh@math.kth.se}
\urladdr{http://www.math.kth.se/~henriksh/}
\author[G.S. Weiss]{Georg S. Weiss}
\address{Graduate School of Mathematical Sciences,
University of Tokyo, 3-8-1 Komaba, Meguro-ku, Tokyo-to, 153-8914 Japan,}
\email{gw@ms.u-tokyo.ac.jp}
\urladdr{http://www.ms.u-tokyo.ac.jp/~gw/}
\thanks{$2000$ {\it Mathematics Subject Classification.\/} Primary
35R35, Secondary 35B40, 35J60.}
\thanks{H. Shahgholian has been supported in part by
the Swedish Research Council.
G.S. Weiss has been partially supported by the Grant-in-Aid
18740086 of the Japanese Ministry of Education, Culture, Sports, Science and Technology.
He also thanks the Knut och Alice Wallenberg foundation for a visiting appointment to
KTH.
Both J. Andersson and G.S. Weiss thank the G\"oran Gustafsson Foundation
for visiting appointments to
KTH.
The present result is part of the ESF-program GLOBAL. It was completed while the first two authors were visiting the Petrolium Institute in Abu Dhabi.}
\date{}

\theoremstyle{plain}
\newtheorem{theorem}{Theorem}[section]

\newtheorem*{theoremb}{Theorem A}

\newtheorem{lemma}[theorem]{Lemma}
\newtheorem{proposition}[theorem]{Proposition}
\newtheorem{corollary}[theorem]{Corollary}
\theoremstyle{definition}
\newtheorem{definition}[theorem]{Definition}
\theoremstyle{example}

\theoremstyle{definition}
\newtheorem{remark}[theorem]{Remark}
\numberwithin{equation}{section}

\def\R{{\mathbb R}}

\def\N{{\mathbb N}}

\def\dist{\textrm{\rm dist}}

\def\dh{d{\mathcal H}^{n-1}}
\def\dho{d{\mathcal H}^1}

\begin{document}
\maketitle
\begin{abstract}
We introduce a new method for the analysis of singularities
in the unstable problem
$$\Delta u = -\chi_{\{u>0\}},$$
which arises in solid combustion as well as in the composite membrane problem.
Our study is confined to points of ``{\em supercharacteristic}'' growth of the solution, i.e.
points at which the solution grows faster than the characteristic/invariant
scaling of the equation would suggest.
At such points the classical theory is doomed to fail, due to
incompatibility of the invariant scaling of the equation and the
scaling of the solution.
\\
In the case of
two dimensions
our result shows that
in a neighborhood of the
set at which the second derivatives of $u$ are unbounded,
the level set $\{ u=0\}$ consists of two $C^1$-curves meeting at right angles.
It is important that our result is not confined to the minimal solution
of the equation but holds for {\em all} solutions.
\end{abstract}
\tableofcontents
\section{Introduction}
In the last decade, the theory of free boundary regularity of obstacle type has got renewed attention, owing to the
seminal paper \cite{revis} of L.A. Caffarelli as well as
\cite{cks}.
Many interesting old and new problems,
intractable by earlier techniques,
have been solved, thanks to the ideas in \cite{revis} and \cite{cks} (see for example \cite{suw}).
All these problems share a common feature: the scaling of the solution at
free boundary points {\em coincides} with the characteristic/invariant scaling of the equation.
However, there are problems arising in applications for which this does
not hold. An example is the unstable obstacle problem
\begin{equation}\label{eq}
\Delta u = -\chi_{\{ u>0\}}\;\quad  \hbox{ in } \Omega\subset \R^n,
\end{equation}
related to traveling wave solutions in solid
combustion with ignition temperature (see the introduction of
\cite{monneauweiss} for more details), to
the composite membrane problem (see \cite{chanillo1},
\cite{chanillo2}, \cite{blank}, \cite{shahgholian}, \cite{chanillokenig}, \cite{chanillokenigtu}) as well as
the shape of self-gravitating rotating
fluids describing stars (see \cite[equation (1.26)]{stars}).
Solutions of equation (\ref{eq}) may exhibit ``{\em supercharacteristic}''
growth of order $$r^2 |\log r|$$ not suggested by the invariant/characteristic scaling $u(rx)/r^2$ of the equation.\\
In this paper we introduce a new method
to analyze the fine structure of singular sets close to points
of supercharacteristic growth of the solution.
\\
Equation (\ref{eq}) has been investigated by R. Monneau-G.S. Weiss
in \cite{monneauweiss}. They establish partial regularity for
\textsl{second order non-degenerate} solutions of (\ref{eq}).
More precisely they show that the singular set has
Hausdorff dimension less than or equal to $n-2$, and that in two
dimensions the free boundary consists close to points where the
second derivative is unbounded, of four Lip\-schitz graphs meeting
at right angles. They also show that energy-minimising solutions are
in the two-dimensional case of class $C^{1,1}$
and that their free boundaries are locally analytic.\\
J. Andersson-G.S. Weiss have constructed a cross-shaped counter-example
proving that the solution need not be of class $C^{1,1}$ (see \cite{anderssonweiss}).
In \cite{monneauweiss} it has been shown that the second variation
of the energy at that particular solution takes the value $-\infty$.
In this sense the cross-solution
is completely unstable. Moreover, it cannot be obtained by naive
numerical schemes.
\\
In this paper we analyze the behavior of solutions
at points at which the second derivatives are unbounded.
Difficulties in the analysis are:\\
(i) At cross-like singular points the solution has the ``wrong scaling'', i.e.
$u(rx)$ scales like $r^2 |\log(r)|$ which is different from the
characteristic scaling $r^2$ of the equation.
The lack of a suitable local Lyapunov functional/monotonicity formula
implies that methods
like the Lojasiewicz inequality (see for example \cite{simon}, \cite{simon2}) would be hard to apply
even at isolated singularities.
\\
(ii) The cross-like singularities are unstable.
\\
(iii) The comparison principle does not hold.
\\
Instead we use knowledge about the Newtonian potential
of the right-hand side to derive a
quantitative estimate for the projection
of the solution onto the homogeneous
harmonic polynomials of degree $2$.
This leads in the case of two dimensions to 
the growth estimate Theorem A (i) for the solution
as well as
an estimate
of order
\begin{equation}\label{letscallitstar}
\int_0^r \frac{\sqrt{|\log|\log s||}}{s|\log s|^{3/2}}\> ds
\end{equation}
for how much the projection of $u(x+s\cdot)$ and also the
approximate tangent space of the singular set can turn as $s$ moves
from $r$ to $0$ (see Theorem A and Remark \ref{intrem}). Our main result Theorem A shows
that close to a non-degenerate singular point, the level set $\{
u=0\}$ consists of two $C^1$-curves meeting at right angles. We
provide estimates for the modulus of the normal of the free boundary
close to singular points. 
Different from the (also two-dimensional) unique tangent cone result 
\cite[Theorem 7.1]{monneauweiss},
the result in the present paper is a {\em quantitative} result
valid {\em uniformly} for a certain class of solutions.
Moreover the result in the present paper is not confined to
the minimal solution. 
\\
In the paper \cite{lines3d} in preparation the authors
extend these new methods to the case of higher dimensions.
\\
Our main result in the present paper is the following (cf. Corollary \ref{logg} and Corollary \ref{iiandiii}):
\begin{theoremb}\label{mainthm}
Let $u$ be a solution of (\ref{eq}) in $\Omega\subset \R^2$
satisfying $\sup_{\Omega} |u|\le M$. Moreover let $d>0$. Then there
exist an $r_0=r_0(M,d)>0$ and a $\delta_0=\delta_0(M,d)>0$ such that
if $x^0\in \Omega_d= \{x\in \Omega:\; \textrm{dist}(x,\partial
\Omega)>d\}$ and
\begin{equation}\label{letscallit3star}
S^u(x^0,r)\equiv \Big( \frac{1}{r^{n-1}}\int_{\partial B_r(x^0)} u^2\dho
\Big)^{1/2}\ge \frac{r^2}{\delta}
\end{equation}
for some $\delta\le \delta_0$, $r\le r_0$ and $u(x^0)=|\nabla
u(x^0)|=0$ then:

(i) $\big(\frac{1}{\delta}-C(M,d)\big)s^2+c\log(r/s)s^2\le S^u(x^0,s)$ for
every $s\le r$.

(ii) There exists a second order homogeneous harmonic polynomial
$p^{x^0,u}=p$ such that for each $\alpha\in (0,1/2)$ and each $\beta\in (0,1)$,
\begin{equation}\label{letscallit2star}
\Big{\|} \frac{u(x^0+sx)}{\sup_{B_s(x^0)}|u|}-p
\Big{\|}_{C^{1,\beta}}\le C(M,d,\alpha,\beta) \Big(
\frac{\delta}{1+\delta\log(r/s)}\Big)^{\alpha}.
\end{equation}

(iii) The set $\{u=0\}\cap B_r(x^0)$ consists of two $C^1$-curves intersecting
each other at right angles at $x^0$.
\end{theoremb}
\begin{remark}\label{intrem}
1) By \cite[Lemma 8.5]{monneauweiss} the estimate Theorem A (i) is sharp.
The inequality (\ref{letscallit3star}) is always satisfied for some
$r$ at singular points, that is, points at which the solution $u$ is
not $C^{1,1}$. Theorem A thus states that $x^0$ is a singular point
if and only if (\ref{letscallit3star}) is satisfied for some $r$.

2) The left hand side in
(\ref{letscallit2star}) may be estimated by the somewhat sharper term in
(\ref{letscallitstar}) (see the end of the proof of Theorem \ref{auto}).
\end{remark}

The proof of {\sl (i)} in Theorem A is contained in Corollary
\ref{logg}, and {\sl (ii)} and {\sl (iii)} will be proved in
Corollary \ref{iiandiii}.

\section{Notation}
Throughout this article $\R^n$ will be equipped with the Euclidean
inner product $x\cdot y$ and the induced norm $\vert x \vert\> .$
We define $e_i$ as the $i$-th unit vector in $\R^n\> ,$ and
$B_r(x^0)$ will denote the open $n$-dimensional ball of center
$x^0\> ,$ radius $r$ and volume $r^n\> \omega_n\> .$
When not specified, $x^0$ is assumed to be $0$.
We shall often use
abbreviations for inverse images like $\{u>0\} :=
\{x\in \Omega\> : \> u(x)>0\}\> , \> \{x_n>0\} :=
\{x \in \R^n \> : \> x_n > 0\}$ etc.
and occasionally
we shall employ the decomposition $x=(x_1,\dots,x_n)$ of a vector $x\in \R^n\> .$
Since we are concerned with local regularity we will
use the set $\Omega_d := \{ x \in \Omega: \dist(x,\partial\Omega)\ge d >0\}.$
We will use the $k$-dimensional Hausdorff measure
${\mathcal H}^k$.
When considering a set $A\> ,$ $\chi_A$ shall stand for
the characteristic function of $A\> ,$
while
$\nu$ shall typically denote the outward
normal to a given boundary.
\section{Preliminaries}\label{prelim}
In this section we state some of the definitions and tools
from \cite{cpde}, \cite{monneauweiss} and mention some examples from
\cite{anderssonweiss}.\\

First we need the monotonicity formula derived
in \cite{cpde} by G.S. Weiss for a class of semilinear free boundary
problems. For the sake of completeness let us state
the unstable case here:
\begin{theorem}[Monotonicity formula, \cite{cpde}]
\label{mon}
Suppose that $u$ is a solution of (\ref{eq}) in $\Omega$
and that $B_\delta(x^0)\subset \Omega\> .$
Then for all $0<\rho<\sigma<\delta$
the function
\[ \Phi^u_{x^0}(r) := r^{-n-2} \int_{B_r(x^0)} \left(
{\vert \nabla u \vert}^2 \> -\> 2\max(u,0)
\right)\]\[
- \; 2 \> r^{-n-3}\>  \int_{\partial B_r(x^0)}
u^2 \> d{\mathcal H}^{n-1}\; ,\]
defined in $(0,\delta)\> ,$ satisfies the monotonicity formula
\[ \Phi^u_{x^0}(\sigma)\> -\> \Phi^u_{x^0}(\rho) \; = \;
\int_\rho^\sigma r^{-n-2}\;
\int_{\partial B_r(x^0)} 2 \left(\nabla u \cdot \nu - 2 \>
{u \over r}\right)^2 \; d{\mathcal H}^{n-1} \> dr \; \ge 0 \; \; .\]
\end{theorem}
The following proposition has been proved in \cite[Section 5]{monneauweiss}.
\begin{proposition}[Classification of blow-up limits with fixed center, Proposition 5.1 in \cite{monneauweiss}]\label{fixedcenter}
Let $u$ be a solution of (\ref{eq}) in $\Omega$
and let us consider a point
$x^0\in \Omega\cap \{ u=0\}\cap\{ \nabla u =0\}.$\\
(i) In the case $\Phi^u_{x^0}(0+)=-\infty$,
$\lim_{r\to 0} r^{-3-n}\int_{\partial B_r(x^0)} u^2 \> d{\mathcal H}^{n-1}
= +\infty$, and for
$S^u(x^0,r) = \left(r^{1-n}\int_{\partial B_{r}(x^0)} u^2\> d{\mathcal H}^{n-1}
\right)^{1\over 2},$
each limit of
\[ \frac{u(x^0+r x)}{S^u(x^0,r)}\]
as $r\to 0$ is a homogeneous harmonic polynomial of degree $2$.
\\
(ii) In the case $\Phi^u_{x^0}(0+)\in (-\infty,0)$,
\[ u_r(x) := \frac{u(x^0+r x)}{r^2}\]
is bounded in $W^{1,2}(B_1(0))$,
and each limit as $r\to 0$ is a homogeneous solution of degree $2$.\\
(iii) Else $\Phi^u_{x^0}(0+)=0$, and
\[ \frac{u(x^0+r x)}{r^2}\to 0\hbox{ in } W^{1,2}(B_1(0)) \hbox{ as } r\to 0\; .\]
\end{proposition}
\begin{remark}\label{int}
1. As observed recently by one of the authors, case (ii) is
possible even in two dimensions (cf. \cite{lines3d}).\\
2. Case (iii) is equivalent to $u$ being degenerate of second order
at $x^0$.
\end{remark}
In \cite{anderssonweiss}, the authors have obtained
abstract existence of solutions
in two dimensions that exhibit {\em cross-like singularities},
at which the second derivatives of the solution are
unbounded (case (i) of Proposition \ref{fixedcenter}), as well as degenerate singularities, at which
the solution decays to zero faster than any quadratic
polynomial (case (iii) of Proposition \ref{fixedcenter}):
\begin{theorem}[Cross-shaped singularity in two dimensions, Corollary 4.2 in \cite{anderssonweiss}]\label{cross}
There exists a solution $u$ of
\begin{displaymath}
\Delta u= -\chi_{\{u>0\}} \quad \textrm{in } B_1\subset \R^2
\end{displaymath}
that is {\bf not} of class $C^{1,1}$.
Each limit of
\[ \frac{u(r x)}{S^u(0,r)}\]
as $r\to 0$ coincides after rotation with the function
$(x_1^2-x_2^2)/\Vert x_1^2-x_2^2\Vert_{L^2(\partial B_1(0))}$.
\end{theorem}
\begin{theorem}[Existence of a degenerate point, Corollary 4.4 in \cite{anderssonweiss}]\label{ast}
There exists
a non-trivial solution $u$ of
\begin{displaymath}
\Delta u= -\chi_{\{u>0\}} \quad \textrm{in } B_1\subset \R^2
\end{displaymath}
that is degenerate
of second order at the origin.
\end{theorem}

\section{A Newtonian potential and its projection}\label{newtandproj}
In what follows we will need the space $P$
of second order homogeneous harmonic polynomials and two dimensional
homogeneous polynomials respectively which we define now.
\begin{definition}
Let us first define in each dimension $n\ge 2$ the space
$P$ of $2$-homogeneous harmonic polynomials, i.e. harmonic
polynomials of degree $2$.
\end{definition}

\begin{definition}\label{projection}

(i) Let us define the projection
$$\Pi:W^{2,2}(B_1)\to P$$ as follows: for $v\in W^{2,2}(B_1)$, let
$\Pi(v)$ be the, by Lemma \ref{proj} unique, minimizer of
$$
p\mapsto \int_{B_1}|D^2 v-D^2 p|^2
$$
on $P$,
where $|A| = \sqrt{\sum_{i,j=1}^n a_{ij}^2}$ is the
Frobenius norm of the matrix $A$.
\\
(ii) Let us also define $\tau(v)\ge 0$ by
$$\Pi(v)=\tau(v)p,\; p\in P, \;\sup_{B_1}|p|=1.$$\end{definition}
\begin{lemma}\label{proj}

(i) For each $v\in W^{2,2}(B_1)$ the minimizer of
Definition \ref{projection} exists and is unique.
Thus $\Pi:W^{2,2}(B_1)\to P$ is well-defined.
\\
(ii) $\Pi$ is a linear operator.\\
(iii) If $h\in W^{2,2}(B_1)$ is harmonic in $B_1$
then $\Pi(h(x))=\Pi(h(rx)/r^2)$ for all $r\in (0,1)$.\\
(iv) For every $v,w\in W^{2,2}(B_1)$,
$$\sup_{B_1}|\Pi(v+w)|\le \sup_{B_1}|\Pi(v)|\> +\> \sup_{B_1}|\Pi(w)|.$$
\end{lemma}
\proof The first and second statement follow from the projection
theorem with respect to the $L^2(B_1;\R^{n^2})$-inner product and
the linear subspace $$\{ f\in  L^2(B_1;\R^{n^2})\> : \> f(x)
\textrm{is symmetric, constant, and } \textrm{trace}(f) = 0\}.$$
\\
Writing $h$ as the sum of homogeneous harmonic polynomials $h_j$
that are orthogonal to each other with respect to
$$(v,w) := \int_{B_1} \sum_{i,j=1}^n \partial_{ij} v \partial_{ij} w,$$
we see that $\Pi(h_j)=0$ for all $j$ such that the degree of $h_j$
is different from $2$, implying the third statement.
\\
The last statement follows from the linearity of $\Pi$ and
the triangle inequality in $L^2(B_1;\R^{n^2})$.
\qed\\
In \cite{karpmargulis} L. Karp-A.S. Margulis derive eigenfunction expansions for
generalized Newtonian potentials with respect to a large class of right-hand
sides. In the following lemma we calculate explicitly
a normalized generalized Newtonian potential of $-\chi_{\{ x_1x_2>0\}}$
as well as its projections. Properties (iv), (v) and (vi) in Lemma \ref{limsol}
are crucial for what follows.
\begin{lemma}\label{limsol}
Define $v:(0,+\infty)\times [0,+\infty)\to \R$ by
$$
v(x_1,x_2):=
-4x_1x_2\log(x_1^2+x_2^2)+2(x_1^2-x_2^2)\left( \frac{\pi}{2}-2\arctan\left(\frac{x_2}{x_1}\right) \right)-\pi(x_1^2+x_2^2).
$$
Moreover let
$$
w(x_1,x_2):=\left\{ \begin{array}{ll}
v(x_1,x_2),& x_1x_2\ge 0, x_1\ne 0,\\
-v(-x_1,x_2),& x_1<0, x_2\ge 0,\\
-v(x_1,-x_2),& x_1>0, x_2\le 0,
\end{array}\right.$$
and let
$$ z(x_1,x_2) := \frac{w(x_1,x_2)-\pi(x_1^2+x_2^2)+8x_1x_2}{8\pi}.
$$
Then, $z$ is the unique solution to\\
(i) $\Delta z=-\chi_{\{x_1x_2>0\}}$ in $\R^2$,\\
(ii) $z(0)=|\nabla z(0)|=0$,\\
(iii) $\lim_{x \to \infty}{z(x)\over {|x|^3}}=0$,\\
(iv) $\Pi(z)=0$,\\
(v) $\Pi(z_{1/2})=\log(2)x_1x_2/\pi,$\\
(vi) $\tau(z_{1/2})=\log(2)/(2\pi)$.
\end{lemma}
\proof
A calculation shows that $w$ can be extended to a $C^1$-function
and that $\Delta w = -4\pi\chi_{\{x_1x_2>0\}} + 4\pi\chi_{\{x_1x_2<0\}}$.
We obtain that $z$ can be extended to a $C^1$-function
solving $\Delta z = -\chi_{\{x_1x_2>0\}}$ in $\R^2$ and satisfying (ii) and (iii).\\
Next we show that $h:=\Pi(z)=0$:
setting $$D^2 h = \begin{pmatrix}a&b\\
b&-a\end{pmatrix},$$
we obtain
$$0=\partial_b \int_{B_1}|D^2 z-D^2 h|^2
= 4 \int_{B_1} \partial_{12}(h-z)
= 4b - 4 \int_{B_1} \partial_{12} z
$$ $$= 4b + 2 \int_{B_1} {1+\log (x_1^2+x_2^2)\over \pi}
= 4b$$
as well as
$$ 0=\partial_a \int_{B_1}|D^2 z-D^2 h|^2
= 4a,$$
implying that $h\equiv 0$.
\\
Rescaling $z$ we see that
$${z(rx_1,rx_2)\over {r^2}}
= z(x_1,x_2) - {x_1x_2 \log r^2  \over {2\pi}}$$
which implies
$$\Pi(z_{1/2})
=\Pi(z) - \Pi({x_1x_2 \log(\left({1\over 2}\right)^2)  \over
{2\pi}}) = -\log(1/2)\Pi(x_1x_2)/\pi = -\log(1/2)x_1x_2/\pi.$$ Thus
{\sl (v)} and {\sl (vi)} are true.

Last, we show uniqueness of $z$ satisfying {\sl (i)-(iv)}. Observe
that {\sl (v)} and {\sl (vi)} are not needed to show uniqueness. If
$z^1$ and $z^2$ are two solutions to {\sl (i)-(iv)}, then by {\sl
(i)}, $z^1-z^2$ is harmonic. Condition {\sl (iii)} implies that
$z^1-z^2$ is a second order polynomial. Conditions {\sl (ii)}
and {\sl (iv)} then imply that $z^1-z^2=0$.
\qed

\section{Growth of the Solution at Singular Points.}

The next lemma is crucial for all that follows.

\begin{lemma}\label{normblow} Let $u$ solve (\ref{eq}) and suppose that $d>0$,
$\sup_\Omega |u|\le M<+\infty$, $x^0\in \Omega_d$, $u(x^0)=|\nabla u(x^0)|=0$
and $r\le d/2$. Then
$$
\Big(\int_{B_1}\Big| D^2\frac{u(x^0+rx)}{r^2}-D^2\Pi\big(\frac{u(x^0+rx)}{r^2} \big)\Big|^p \Big)^{1/p}\le C(n,M,d,p)
$$
and
$$
\Big{\|} \frac{u(x^0+rx)}{r^2}-\Pi\big(\frac{u(x^0+rx)}{r^2} \big)\Big{\|}_{C^{1,\beta}}\le C(n,M,d,\beta).
$$
\end{lemma}
\proof Let $u_r(x)=\frac{u(x^0+rx)}{r^2}$. From \cite[4.1
Proposition 1]{stein} we infer that $D^2 u$ is locally of class BMO,
and that
$$
\Big(\int_{B_{3/2}}|D^2 u_r- \overline{D^2 u_{3r/2}}|^2\Big)^{1/2} \le C_1,
$$
where
$$
\overline{D^2 u_{3r/2}}=\frac{1}{\omega_n {(3/2)}^n}\int_{B_{3/2}} D^2 u_r,
$$
and $C_1$ is a constant depending only on $n$, $M$ and $d$. It follows that
$$
C_1\ge \Big( \int_{B_{3/2}}|D^2 u_r- \overline{D^2
u_{3r/2}}|^2\Big)^{1/2}
$$
$$
\ge \Big( \int_{B_{3/2}}|D^2 u_r-
\big(\overline{D^2 u_{3r/2}}-\frac{1}{n}\textrm{trace}(\overline{D^2 u_{3r/2}})I\big)|^2\Big)^{1/2}
$$
$$
-\Big( \int_{B_{3/2}}|
\frac{1}{n}\textrm{trace}(\overline{D^2 u_{3r/2}})I|^2\Big)^{1/2},
$$
where $I$ is the identity matrix. Next it is easy to see that
$$
\int_{B_{3/2}}|
\frac{1}{n}\textrm{trace}(\overline{D^2 u_{3r/2}})I|^2\le 1,
$$
since
$$
\textrm{trace}\big(\overline{D^2 u_{3r/2}}\big)
=\frac{1}{\omega_n(3/2)^n}\int_{B_{3/2}} \Delta u_r
$$
and $|\Delta u_r|\le 1$. In particular we have
$$
C_1+1\ge \Big( \int_{B_{3/2}}|D^2 u_\eta-
\big(\overline{D^2 u_{3r/2}}-\frac{1}{n}\textrm{trace}(\overline{D^2 u_{3r/2}})I\big)|^2\Big)^{1/2}.
$$
Using the minimizing property of the projection $\Pi$ we get
$$
(C_1+1)^2\ge \int_{B_{3/2}}|D^2 u_r-
\big(\overline{D^2 u_{3r/2}}-\frac{1}{n}\textrm{trace}(\overline{D^2 u_{3r/2}})I\big)|^2
$$
$$
\ge \int_{B_{3/2}}|D^2 u_r-D^2 \Pi(u_{3r/2})|^2.
$$
Observe that if we set $v := u_r-\Pi(u_{3r/2})$, then
$$\int_{B_{3/2}}|D^2 v|^2 \le (C_1+1)^2,
 \Vert \Pi(v)\Vert_{L^2(B_1)}\le C_2,$$
$$ \Vert \Pi(v)\Vert_{L^2(B_{3/2})}\le C_3
\textrm{ and } \Vert v-\Pi(v)\Vert_{L^2(B_{3/2})}\le C_4.$$
It follows that $D^2\big( u_{r}-\Pi(u_{r})\big)$ is bounded in $L^2(B_{3/2})$.
Moreover, since $\Pi(u_r)$ is harmonic,
$\Delta \big( u_r-\Pi(u_r) \big)=-\chi_{\{u_r>0\}}$. Poincare's inequality
implies that
$$
\big{\|} u_r-\Pi(u_r)-\overline{\nabla u_r}\cdot x-
\overline{u_r} \big{\|}_{W^{2,2}(B_{3/2})}\le \big{\|}D^2 u_r - D^2\Pi(u_r)\big{\|}_{L^2(B_{3/2})}\le C_5,
$$
where $\overline{\nabla u_r}$ and $\overline{u_r}$ denote the averages. Thus
$L^p$-theory (see for example \cite[Theorem 9.11]{gt}) implies that
$$
\big{\|} u_r-\Pi(u_r)-\overline{\nabla u_r}\cdot x- \overline{u_r} \big{\|}_{W^{2,p}(B_{1})}\le 
C_6.$$
The embedding into H\"older spaces therefore yields
$$
\big{\|} u_r-\Pi(u_r)-\overline{\nabla u_r}\cdot x- \overline{u_r} \big{\|}_{C^{1,\beta}(B_{1})}\le C_7.
$$
Using that $u(x^0)=|\nabla u(x^0)|=0$ and the above estimates implies the
statement of the Lemma.\qed

\begin{remark}\label{comparable}
The above Lemma implies in particular that when one of the quantities
$\|u\|_{L^\infty(B_r(x^0))}$, $S^u(x^0,r)$ and $\tau(u(x^0+r\cdot))$
is large in comparison to $r^2$ then all these quantities are
comparable. Let us indicate how to prove this: assume that
$\tau(u(x^0+r\cdot))> \bar{C} r^2$ for some large constant $\bar{C}=\bar{C}(n,M,d)$ then
$$
S^u(x^0,r)=\Big(\frac{1}{r^{n-1}}\int_{\partial B_r(x^0)} u^2\dh
\Big)^{1/2}\ge
 \Big(\frac{1}{r^{n-1}}\int_{\partial B_r(x^0)}
\Pi(u)^2 \dh\Big)^{1/2}
$$
$$
-\Big(\frac{1}{r^{n-1}}\int_{\partial B_r(x^0)} (u-\Pi(u))^2\dh
\Big)^{1/2}\ge c(n)\tau(u(x^0+r\cdot))-C(n,M,d) r^2.
$$
It follows that if $\bar{C} >  2C(n,M,d)/c(n)$ then $S^u(x^0,r)>
c(n)\tau(u(x^0+r\cdot))/2.$ Similarly one may deduce that 
under the above assumptions
$S^u(x^0,r)<
C(n)\tau(u(x^0+r\cdot))$ and that the corresponding relationships between the
other quantities above hold.
\end{remark}

In what follows, we denote by $z(x_1,\dots,x_n):= z(x_1,x_2)$ the solution of Lemma \ref{limsol}, extended to $\R^n$.
\begin{lemma}\label{zconv}
For each $\epsilon>0,n\in \N, d>0,M<+\infty,\alpha \in [1,+\infty)$ and $\beta\in (0,1)$ there exist $r_0,\delta>0$ with the following property:\\
Suppose that $0<r\le r_0$, $x\in\Omega_d$ and that
$u$ is a solution of (\ref{eq}) in $\Omega$ satisfying $\sup_{\Omega} |u|\le M$, $u(x)=|\nabla u(x)|=0$ and
$${\mathcal L}^n((\{ u(x+r\cdot)>0\} \triangle \{ x_1x_2>0\})\cap B_1)\le \delta.$$
Then
$$
\left\Vert{u(x+{r}\cdot)\over {{r}^2}}-\Pi({u(x+{r}\cdot)\over {{r}^2}})
- z\right\Vert_{C^{1,\beta}(\bar B_1)} \le \epsilon.$$
\end{lemma}
\proof
Suppose that $r_j\to 0$, that
$${\mathcal L}^n(\{ u_j(x^j+r_j\cdot)>0\} \triangle \{ x_1x_2>0\})\to 0\textrm{ as } j\to\infty$$
and that
$$
{u_j(x^j+{r_j}\cdot)\over {{r_j}^2}}-\Pi({u_j(x^j+{r_j}\cdot)\over {{r_j}^2}})
\to \tilde z
\textrm{ in }C^{1,\beta}_{loc}(\R^n)
\textrm{ and weakly in } W^{2,\alpha}_{loc}(\R^n)$$
as $j\to\infty$
(cf. Lemma \ref{normblow}).
\\
Now let $\tilde N$ be the Newtonian potential of $\chi_{\Omega_d}\Delta u_j$,
i.e.
$$\tilde N(y) :=
\left\{\begin{array}{ll}
{1\over {n(2-n)\omega_n}} \int_{\R^n} |y-\xi|^{2-n} (\chi_{\Omega_d}\Delta u_j)(\xi)\> d\xi,&n>2,\\
{1\over {2\pi}} \int_{\R^2} \log |y-\xi| (\chi_{\Omega_d}\Delta u_j)(\xi)\> d\xi,&n=2.
\end{array}\right.$$
Next we let $N(y) := \tilde N(y)-\tilde N(x^j)-\nabla \tilde N(x^j)\cdot (y-x^j)$, and
consider the harmonic function $h(y):= u_j(y)-N(y)$.
Since $\sup_{\Omega} |u_j|\le M$, $|h|\le C_2$ on $\partial B_d(x^j)$,
and it follows that $|D^3 h(y)| \le C_3$ in $B_{d/2}(x^j)$, where $C_3$ depends
on $n$, $d$ and $M$.
Consequently
$$|u_j(y)-N(y)-D^2 h(x^j)(y-x^j)(y-x^j)|\le C_4 |x^j-y|^3\textrm{ in }B_{d/2}(x^j),$$
where $C_4$ depends only on $n,d$ and $M$.
For the scaled functions $v_j(y):=u_j(x^j+r_jy)/r_j^2, N_j(y):=N(x^j+r_jy)/r_j^2$
and $p_j(y) = D^2 h(x^j)(y)(y)$ we obtain
$$|v_j(y)-N_j(y)-p_j(y)|\le C_4 r_j |y|^3\textrm{ in }
B_{d/(2r_j)}.$$ Thus
$$v_j-\Pi(v_j) = N_j - \Pi(N_j) + o(1)\textrm{ as } j\to\infty.$$
Passing if necessary to another subsequence $j\to\infty$,
the functions $N_j$ converge locally to $N_0$, where
$$
\Delta N_0=-\chi_{\{ x_1x_2>0\}}, N_0(0)=0, \nabla N_0(0)=0 \textrm{ and } N_0-\Pi(N_0)=\tilde z.
$$
We need to establish that $|N_0(y)|=o(|y|^3)$ as $|y|\to\infty$.
Once this is established the uniqueness part of Lemma \ref{limsol}
implies that $\tilde z = N_0-\Pi(N_0)=z$ and the Lemma follows.
First, $D^2 N_0 \in BMO(\R^n)$, so that
$$\int_{B_1} \left|
\frac{D^2(N_0(Ry))-\overline{D^2(N_0(R\cdot))}}{\sup_{B_R}|D^2N_0|}\right|^2\>
dy \le C_5 \frac{R^4}{\sup_{B_R}|D^2N_0|^2}$$ for all $R\in
(0,+\infty)$, where $\overline{D^2(N_0(R\cdot))}$ denotes the mean
value of $D^2(N_0(R\cdot))$ on $B_1$. Thus
$\limsup_{R\to\infty}\sup_{B_1}|D^2N_0(R\cdot)|/R^2= +\infty$
implies that
\begin{equation}\label{N0}
\begin{array}{l}
N_0(R_k\cdot)/\sup_{B_{R_k}}|D^2N_0| \textrm{ converges for a sequence }
R_k\to\infty\\ \textrm{ to a }2 \textrm{-homogeneous harmonic polynomial.}
\end{array}
\end{equation}
Now suppose towards a contradiction that
$$\limsup_{|y|\to\infty} \frac{|N_0(y)|}{|y|^3}>0.$$
Then $\Delta (N_0-z)=0$ in $\R^n$ and
$$\limsup_{|y|\to\infty} \frac{|N_0(y)-z(y)|}{|y|^3}>0.$$
Thus $N_0-z$ must be a harmonic polynomial of degree $m\ge 3$,
contradicting (\ref{N0}).
\qed

\begin{lemma}\label{omegdiff}
Let
$n=2, d>0$ and $M<+\infty$. Then there are $r_0, \delta>0$ with the following property:\\
Suppose that $0<r\le r_0,x^0\in\Omega_d$ and that $u$ is a solution
of (\ref{eq}) in $\Omega$ satisfying $\sup_{\Omega} |u|\le M$,
$u(x^0)=|\nabla u(x^0)|=0$ and
$$
S^u(x^0,r)\ge \frac{r^2}{\delta},
$$
for some $r\le r_0$. Then
$${\mathcal L}^n\big((\{u(x^0+r\cdot) > 0 \} \Delta \{\Pi(u(x^0+r\cdot))>0\})\cap B_1\big)\le C \frac{|\log(S^u(x^0,r)/r^2)|}{S^u(x^0,r)/r^2},
$$
where $C=C(d,M,r_0)$.
\end{lemma}
\proof Let $u_r(y) := u(x^0+ry)/r^2$. Then $u_r$ is a solution to
(\ref{eq}) and $S^{u_r}(0,1)>1/\delta$. Let $\tau(u_r) p_r =
\Pi(u_r)$. By Lemma \ref{normblow}, $\sup_{B_1} |u_r-\tau(u_r)
p_r|\le C$, and we obtain at each point $x \in \{ u_r>0\}\cap \{ p_r\le
0\}$ that
$$
|p_r(x)| \le\frac{C}{\tau(u_r)}\le  \frac{C_1}{S^{u_r}(0,1)},
$$
where we have used that $S^{u_r}(0,1)$ is comparable to $\tau(u_r)$
(see Remark \ref{comparable}). Next we calculate
$${\mathcal L}^n(\{ u_r>0\}\cap \{ p_r\le 0\}\cap B_1)
\le{\mathcal L}^n(\{|p_r|\le \frac{C_1}{S^{u_r}(0,1)}\}\cap B_1)$$
$$
\le 4 {\mathcal L}^n(\{ (x_1,x_2):0<x_1<1, 0<x_2<1, x_1x_2\le
\frac{C_1}{S_r(0,1)}\})$$ $$ = 4 \int_0^{C_1/S^{u_r}(0,1)} \> dx_1 +
4 \int_{C_1/S^{u_r}(0,1)}^1 {C_1\over {x_1S^{u_r}(0,1)}} \> dx_1 \le
\frac{C|\log(S^{u_r}(0,1))|}{S^{u_r}(0,1)}.$$ The Lemma follows by
scaling back $S^u(x^0,r)=r^2S^{u_r}(0,1).$ \qed

\begin{lemma}\label{growth}
Let $n=2$. For each $\gamma \in (0,\log(2)/(2\pi)), d>0$ and
$M<+\infty$
there are $r_0,\delta>0$, depending only on $\gamma$, $d$ and $M$,
with the following property:\\
Suppose that $0<r\le r_0,$  $x^0\in\Omega_d$ and that $u$ is a
solution of (\ref{eq}) in $\Omega$ satisfying $\sup_{\Omega} |u|\le
M$, $u(x^0)=|\nabla u(x^0)|=0$ and for some $r\le r_0$,
$$
S^u(x^0,r)\ge \frac{r^2}{\delta}.
$$
Then $\tau(4u(x^0+r\cdot/2)/r^2)\ge \tau(u(x^0+r\cdot)/r^2)+\gamma$.
\end{lemma}
\proof
Suppose towards a contradiction that
$\tau(4u_j(x^j+r_j\cdot/2)/r_j^2)<\tau(u_j(x^j+r_j\cdot)/r_j^2)+\gamma$
for a sequence $u_j$ satisfying the assumptions with $\delta=\delta_j\to 0$ as $j\to\infty$.
Let $v_j := u_j(x^j+r_j\cdot)/r_j^2$.
A straightforward calculation shows that $v_j$ solves (\ref{eq}) and that
$$
S^{v_j}(0,1)\ge \frac{1}{\delta_j}.
$$
From Lemma \ref{omegdiff} it follows that
$$
{\mathcal L}^n(\{v_j > 0 \} \Delta \{\Pi(v_j)>0\})\cap B_1)\to 0.
$$
We may apply Lemma \ref{zconv} and deduce that, after a rotation of
the coordinate system, $v_j -\Pi(v_j)\to z$ weakly in
$W^{2,\alpha}(B_1)$ and strongly in $C^{1,\beta}(\bar B_1)$ as
$j\to\infty$, and that therefore --- rotating each $v_j$ only slightly
more --- $\Pi(v_j)=M_j x_1x_2$ with $M_j\to+\infty$ as $j\to\infty$.
Defining $f_{1/2}(y) := 4f(y/2)$, it follows from Lemma \ref{limsol}
{\sl (v)} that $\Pi((v_j)_{1/2}-M_j x_1x_2)\to
\Pi(z_{1/2})=\log(2)x_1x_2/\pi$ as $j\to\infty$. On the other hand,
$\tau((v_j)_{1/2})< \tau(v_j)+\gamma$, so that
$$(\log(2)/\pi+M_j)/2=\tau((\log(2)/\pi+M_j)x_1x_2)$$ $$=o(1)+\tau((v_j)_{1/2})< o(1)+\tau(v_j)+\gamma
= o(1)+M_j/2 + \gamma,
$$
a contradiction for large $j$. \qed

The next Corollary proves the first statement in Theorem A and is fundamental for the rest of the paper. 

\begin{corollary}\label{logg}
Let $n=2$. Fix a $\gamma\in (0,\log(2)/2\pi)$ and let $u$ satisfy
the assumptions in Lemma \ref{growth} for some $r\le r_0$ (with
possibly somewhat smaller $\delta$). Then
$$
\tau(2^{2j} u(x^0+2^{-j}r\cdot )/r^2)\ge \tau(u(x^0+r\cdot)/r^2)+j\gamma
\textrm{ for all }j\in \mathbb{N}.
$$
Moreover, for each $s\le r$,
$$
\frac{S^u(x^0,s)}{s^2}\ge \frac{S^u(x^0,r)}{r^2}+c\gamma \frac{\log(r/s)}{\log(2)}-2C,
$$
where $c=\|x_1x_2\|_{L^2(\partial B_1)}$, and $C=C(M,d,r_0)$.
\end{corollary}
\proof Since by Lemma \ref{normblow}
$$
\sup_{B_1}\Big| \frac{u(x^0+rx)}{r^2}-\Pi\big(\frac{u(x^0+rx)}{r^2}
\big)\Big|\le C_0,
$$
it follows that for $s\le r$,
$u_s(x)=u(x^0+sx)/s^2$ and 
$c=\|x_1x_2\|_{L^2(\partial B_1)}$,
\begin{align}\label{comptauS1}
S^{u_s}(0,1) - \sqrt{2C_0\pi}
\le & \left( \int_{\partial B_1} |\Pi(u_s)|^2\dho\right)^{1\over 2}
+ \left( \int_{\partial B_1} |u_s-\Pi(u_s)|^2\dho\right)^{1\over 2}\\
\nonumber
&- \sqrt{2C_0\pi}
\le c \tau(u_s).
\end{align}
Similarly it follows that
\begin{equation}\label{comptauS2}
\tau(u_s)\> \Big( \int_{\partial
B_1} (x_1x_2)^2 \Big)^{1/2}\dho\le S^{u_s}(0,1)+\sqrt{2C_0\pi}.
\end{equation}
From Lemma \ref{growth} we infer that if $S^u(x^0,r)/r^2\ge
1/\sigma$ with $\sigma<\delta$ and $\delta$ is as in Lemma
\ref{growth}, then $\tau_{r/2}\ge \tau_r+\gamma$. Here we use
short hand $\tau_r\equiv \tau(u(x+r\cdot)/r^2)$. From inequalities
(\ref{comptauS1}) and (\ref{comptauS2}) we see that
\begin{equation}\label{firstitforgrowth}
\frac{S^u(x^0,r/2)}{(r/2)^2}\ge (\tau_r +\gamma)c-\sqrt{2C_0\pi}\ge
\frac{S^u(x^0,r)}{r^2}+\gamma c -2\sqrt{2C_0\pi},
\end{equation}
where $c$ is the constant in the statement of the Corollary. In
particular, if $\sigma$ has been chosen small enough, say $1/\sigma>
1/\delta+2C_1$, then $u$ satisfies the assumptions of Lemma
\ref{growth} in $B_{r/2}$. We may thus apply Lemma \ref{growth}
again and deduce that
$$
\frac{S^u(x^0,r/4)}{(r/4)^2}\ge (\tau_r +2\gamma)c-2\sqrt{2C_0\pi}.
$$ 
Applying Lemma \ref{growth} $j$ times, we arrive at
$$
\frac{S^u(x^0,r/2^j)}{(r/2^j)^2}\ge (\tau_r +j\gamma)c-C_1\ge
\frac{S^u(x^0,r)}{r^2}+ c\gamma j -2C_1.
$$
Notice that since $\tau_{2^{-j}r}$ is increasing in $j$ and thus
$S^u(x^0,2^{-j}r)\ge \tau_r-2\sqrt{2C_0\pi}$ for each $j$ and 
the assumptions of Lemma \ref{growth} are therefore satisfied for 
each $j$.

If we put $s=2^{-j}r$ then $j=\log(r/s)/\log(2)$ and we obtain the
statement in the Corollary. For general $s\le r$ we may consider a
$j$ such that $2^{-(j+1)}r< s \le 2^{-j}r$. Using Lemma
\ref{normblow},
$$
\Big{\|} \frac{u(x^0+2^{-j}rx)}{(2^{-j}r)^2}-\Pi\big(
\frac{u(x^0+2^{-j}rx)}{(2^{-j}r)^2}
\big)\Big{\|}_{C^{1,\beta}(B_1)}\le C_2,
$$
and it follows that
$$
\Big{|}\frac{S^u(x^0,s)}{s^2}- \frac{S^u(x^0,2^{-j}r)}{(2^{-j}r)^2}
\Big{|}\le C_3.
$$
The Corollary follows with a slightly larger constant $C$.\qed

\section{Controlling the movement of $\Pi(u(x+r\cdot))$}
In this section we will exploit the estimate in Corollary \ref{logg}
to obtain control of how much the projection of $u(x+r\cdot)$
can turn when passing to a smaller radius $r$.

\begin{lemma}\label{gest}
Let $n=2$, $d>0$ and $M<\infty$. Then there is $r_0,\delta>0$ with the following property:\\
Suppose that $0<r\le r_0,$ $x^0\in\Omega_d$ and that $u$ is a solution
of (\ref{eq}) in $\Omega$ satisfying $\sup_{\Omega} |u|\le M$,
$u(x)=|\nabla u(x)|=0$ and
$$
\frac{S^u(x^0,r)}{r^2}\ge \frac{1}{\delta}.
$$
Let $g$ be the solution of
$$
\Delta
g=\chi_{\{\Pi(u(x+r\cdot))>0\}}-\chi_{\{u(x+r\cdot)>0\}}\textrm{ in
}B_1,
$$
$$
g=0 \textrm{ on } \partial B_1.
$$
Then \\
(i) $$\Vert D^2 g\Vert_{L^2(B_1)} \le C
\sqrt{\frac{|\log(S^u(x^0,r)/r^2)|}{S^u(x^0,r)/r^2}}.$$
\\
(ii) 
$$
\tau(g)\le C \sqrt{\frac{|\log(S^u(x^0,r)/r^2)|}{S^u(x^0,r)/r^2}},
$$
where $C=C(d,M,r_0)$.
\end{lemma}
\proof
(i) follows from Lemma \ref{omegdiff} and $L^2$-theory (see for example
\cite[Theorem 8.8]{gt}).\\
(ii) Rotating and setting $p:= \Pi(g) = a_1 x_1^2+a_2x_2^2$,
we obtain
$$
\Vert D^2 p\Vert_{L^2(B_1)} \le C_1 \Vert D^2 g\Vert_{L^2(B_1)}
\le C_2 \sqrt{\frac{|\log(S^u(x^0,r)/r^2)|}{S^u(x^0,r)/r^2}}
$$
and
$$|a_j| \le C_3 \sqrt{\frac{|\log(S^u(x^0,r)/r^2)|}{S^u(x^0,r)/r^2}}$$
for $j=1,2$.\qed\\
The next Proposition already contains the desired estimate
for how much the projection may turn when passing from
$u(x^0+r\cdot)$ to $u(x^0+r\cdot/2)$.

\begin{proposition}\label{rotest}
Let $n=2$, $d>0$ and $M<+\infty$. Then there are $r_0,\delta>0$ with the following property:\\
Suppose that $0<r\le r_0,$ $x^0\in\Omega_d$ and that $u$ is a solution of
(\ref{eq}) in $\Omega$ satisfying $\sup_{\Omega} |u|\le M$,
$u(x)=|\nabla u(x)|=0$ and
$$
\frac{S^u(x^0,r)}{r^2}\ge \frac{1}{\delta}.
$$
Then
$$
\sup_{B_1}\left| \frac{\Pi(u(x+r\cdot))}{\sup_{B_1}|\Pi(u(x+r\cdot))|}-\frac{\Pi(u(x+r\cdot/2))}{\sup_{B_1}|\Pi(u(x+r\cdot/2))|} \right|\le C \frac{\sqrt{|\log(|S^u(x^0,r)/r^2|)|}}{\big( S^u(x^0,r)/r^2\big)^{3/2}},
$$
where $C=C(n,M,d)$.
\end{proposition}
\proof Let us consider $v=u_r-z\circ Q_r-h_r-\tau(u_r)p_r$ where
$u_r(y)=u(x+ry)/r^2$, $\Pi(u_r)=\tau(u_r)p_r$, the orthogonal matrix
$Q_r$ has been chosen such that $\{ \Pi(u_r)>0\}=\{ (x_1x_2)\circ
Q_r>0\}$ (we may assume that $Q_r=I$, the identity matrix), 
$h_r=h(ry)/r^2$, and $h$ is harmonic and satisfies
$h(x)\le C_1|x|^3$. It follows that $\Pi(v)=0$. Moreover we may
express $v=g + \tilde h$ where $g$ is the solution of Lemma
\ref{gest} and $\tilde h$ is harmonic. Lemma \ref{gest} (ii) implies
now that for $\tilde h_{1/2}(y) = 4\tilde h(y/2)$, $g_{1/2}(y) =
4g(y/2)$ and $v_{1/2}(y) = 4v(y/2)$,
$$
\sup_{B_1}|\Pi(v_{1/2})|=\sup_{B_1}|\Pi(\tilde h_{1/2}+g_{1/2})|\le
\sup_{B_1}|\Pi(g_{1/2})|$$ 
$$+\sup_{B_1}|\Pi(\tilde h_{1/2})| \le
\sup_{B_1}|\Pi(\tilde h_{1/2})|+ C_2
\sqrt{\frac{|\log(S^u(x^0,r)/r^2)|}{S^u(x^0,r)/r^2}}.
$$
Since $\Pi(v)=0$ we also know that $|\Pi(\tilde h)|\le |\Pi(g)|\le
C_2\sqrt{\frac{|\log(S^u(x^0,r)/r^2)|}{S^u(x^0,r)/r^2}}$. On the other
hand, using that $\tilde{h}$ is harmonic and Lemma \ref{proj} (iii),
$\Pi(\tilde h)=\Pi(\tilde h_{1/2})$ so that
$$
\sup_{B_1}|\Pi\big( u_{r/2}-z_{1/2}-h_{r/2}-\tau(u_r)p_r\big)|=\sup_{B_1}|\Pi(v_{1/2})|\le
2C_2\sqrt{\frac{|\log(S^u(x^0,r)/r^2)|}{S^u(x^0,r)/r^2}}.
$$
From the linearity of $\Pi$, $|h(x)|\le C_3|x|^3$ and Lemma
\ref{limsol} we infer that
\begin{equation}\label{estPiur}
\sup_{B_1} |\Pi(u_{r/2})-(\tau(u_r)+\log(2)/(2\pi))p_r|
\end{equation}
$$
\le 2C_2
\sqrt{\frac{|\log(S^u(x^0,r)/r^2)|}{S^u(x^0,r)/r^2}}+\sup_{B_1}|\Pi(h_{r/2})|\le
C_4\sqrt{\frac{|\log(S^u(x^0,r)/r^2)|}{S^u(x^0,r)/r^2}};
$$
here we also used that $\sup_{B_1}|\Pi(h_{r/2})|\le C_4r$ which can be
absorbed in the last term since $S^u(x^0,r)/r^2$ is large by assumption. 

From (\ref{estPiur}) we conclude that
$$
\sup_{B_1}\left| \frac{\Pi(u_r)}{\sup_{B_1}|\Pi(u_r)|}-\frac{\Pi(u_{r/2})}{\sup_{B_1}|\Pi(u_{r/2})|} \right|
$$
$$
\le \sup_{B_1}\left| \frac{\Pi(u_r)}{\sup_{B_1}|\Pi(u_r)|}-\frac{(\tau(u_r)+\log(2)/(2\pi))p_r}{\sup_{B_1}|\Pi(u_{r/2})|} \right|+C_6\frac{\sqrt{|\log(S^u(x^0,r)/r^2)|}}{\big(S^u(x^0,r)/r^2\big)^{3/2}},
$$
where we also used $\sup_{B_1}|\Pi(u_{r/2})|\ge C_7 S^u(x^0,r)/r^2$ (c.f. Remark \ref{comparable}). Next we make the following estimate, which together with the
 previous estimate yields the conclusion of the Proposition:
$$
\sup_{B_1}\left| \frac{\tau(u_r) p_r}{\tau(u_r)}-\frac{(\tau(u_r)+\log(2)/(2\pi))p_r}{\sup_{B_1}|\Pi(u_{r/2})|} \right|
$$
$$
\le
\sup_{B_1}\left|
\frac{\tau(u_r)p_r}{\tau(u_r)}-
\frac{(\tau(u_r)+\log(2)/(2\pi))p_r}{(\tau(u_r)+\log(2)/(2\pi))}
\right|
+\left|
\frac{\tau(u_r)+\log(2)/(2\pi)}{\sup_{B_1}|\Pi(u_{r/2})|}-1
\right|
$$
$$
\le C_8\frac{1}{S^u(x^0,r)/r^2}\sqrt{\frac{|\log(S^u(x^0,r)/r^2)|}{S^u(x^0,r)/r^2}},
$$
where we have used (\ref{estPiur}) to estimate
$$|\sup_{B_1} |\Pi(u_{r/2})|-(\tau(u_r)+\log(2)/(2\pi))|\le C_4
\sqrt{\frac{|\log(S^u(x^0,r)/r^2)|}{S^u(x^0,r)/r^2}},
$$
$$
\left|
\frac{\tau(u_r)+\log(2)/(2\pi)}{\sup_{B_1}|\Pi(u_{r/2})|}-1
\right|
$$ 
$$
\le C_9 \frac{1}{S^u(x^0,r)/r^2}\sqrt{\frac{|\log(S^u(x^0,r)/r^2)|}{S^u(x^0,r)/r^2}}.
$$
\qed\\
\begin{theorem}\label{auto}
Let $n=2$, $d>0$ and suppose that $u$ solves (\ref{eq}) and that 
$\sup_{\Omega}|u|\le M<+\infty$. Then
there exists a $\delta=\delta(M,d)>0$ and an
$r_0=r_0(M,d)>0$ such that if $x^0\in \Omega_d$ and
$$
\frac{S^u(x^0,r)}{r^2}\ge \frac{1}{\delta}
$$
for some $r\le r_0$ then for each $\alpha\in(0,1/2)$ and
all $s\le r$,
$$
\sup_{B_1}\Big|\frac{\Pi(u(x^0+rx))}{\sup_{B_1}|\Pi(u(x^0+rx))|}-
\frac{\Pi(u(x^0+sx))}{\sup_{B_1}|\Pi(u(x^0+sx))|} \Big|\le C(d,M,\alpha) \Big(\frac{r^2}{S^u(x^0,r)} \Big)^\alpha.
$$
\end{theorem}
\proof For simplicity we will only prove the Theorem for $s=2^{-j}r$; for
general $s$ we may use the estimate in Lemma \ref{normblow} as indicated in
the proof of Corollary \ref{logg}.

Let us choose $\delta$ small enough so that Corollary \ref{logg}
holds for some fixed $\gamma>0$, i.e.
\begin{equation}\label{twostarin63}
\frac{S^u(x^0,2^{-j}r)}{2^{-2j}r^2}\ge\frac{S^u(x^0,r)}{r^2}+c\gamma j-2C.
\end{equation}
Decreasing $\delta$ somewhat more if necessary, we see that (\ref{twostarin63})
implies that the assumptions in Proposition \ref{rotest} hold for every
ball $B_{2^{-j}r}(x^0)$. Using the triangle inequality we obtain that
$$
\sup_j \bigg[\sup_{B_1}\Big|\frac{\Pi(u(x^0+rx))}{\sup_{B_1}|\Pi(u(x^0+rx))|}-
\frac{\Pi(u(x^0+2^{-j}rx))}{\sup_{B_1}|\Pi(u(x^0+2^{-j}rx))|}
\Big| \bigg]
$$
$$
\le \sum_{j=0}^{\infty}
\bigg[\sup_{B_1}\Big|\frac{\Pi(u(x^0+2^{-j}rx))}{\sup_{B_1}|\Pi(u(x^0+2^{-j}rx))|}-
\frac{\Pi(u(x^0+2^{-j-1}rx))}{\sup_{B_1}|\Pi(u(x^0+2^{-j-1}rx))|}
\Big| \bigg].
$$
This sum may be estimated, by Proposition \ref{rotest},
from above by
\begin{equation}\label{sumthatneedstobeestimated}
\sum_{j=0}^\infty \frac{\sqrt{\log(S^u(x^0,2^{-j}r)/(2^{-2j}r^2))}}{\big( S^u(x^0,2^{-j}r)/(2^{-2j}r^2)\big)^{3/2}}.
\end{equation}
Let us set $k$ to be the smallest integer satisfying
$$
k\ge\frac{1}{c\gamma}\big( \frac{S^u(x^0,r)}{r^2}-2C\big).
$$
For $S^u(x^0,r)/r^2$ large enough we see that \begin{equation}\label{definitionofk}k>c_1 \frac{S^u(x^0,r)}{r^2}.
\end{equation}
Using (\ref{twostarin63}) we may estimate (\ref{sumthatneedstobeestimated})
by
$$
C_2\sum_{j=k}^\infty \frac{\sqrt{\log(c\gamma j)}}{( c\gamma j)^{3/2}}\le C_3 \int_k^\infty
\frac{\sqrt{\log(c\gamma t)}}{( c\gamma t)^{3/2}}dt\le C_4\frac{2+\log{k}}{\sqrt{k}}
\le C_5(\alpha)k^{-\alpha}
$$
for each $\alpha\in (0,1/2)$. Using (\ref{definitionofk}) gives the Theorem.\qed
\section{Conclusion}
\begin{corollary}\label{iiandiii}
Under the assumptions in Theorem \ref{auto} the following holds:

(i) there exists a homogeneous harmonic polynomial $p^{x^0,u}=p$ of second order
such that for each $\alpha\in(0,1/2)$ and each $\beta\in (0,1/2)$
$$
\Big{\|} \frac{u(x^0+sx)}{\sup_{B_s(x^0)}|u|}-p \Big{\|}_{C^{1,\beta}}\le
C(d,M,\alpha,\beta) \Big( \frac{\delta}{1+\delta\log(r/s)}\Big)^{\alpha}.
$$

(ii) The set $\{u=0\}\cap B_r(x^0)$
consists of two $C^1$-curves intersecting each other at right angles at $x^0$.
\end{corollary}
\proof From Corollary \ref{logg} we know that for each $s\le r$
\begin{equation}\label{whateverimaycallit}
\frac{S^u(x^0,s)}{s^2}\ge c_1\big( \frac{1}{\delta}+\log(r/s) \big).
\end{equation}
It follows from Theorem \ref{auto} that
\begin{equation}\label{limmitexists}
\lim_{s\to 0}\frac{\Pi(u(x^0+sx))}{\sup_{B_1}|\Pi(u(x^0+sx))|}=p^{x^0,u}\equiv p
\end{equation}
exists. Using Lemma \ref{normblow} gives
\begin{equation}\label{c1betarotationsetimate}
C_2\ge \Big{\|} \frac{u(x^0+sx)}{s^2}-\frac{\Pi(u(x^0+sx))}{s^2} \Big{\|}_{C^{1,\beta}}
\end{equation}
$$
\ge \big{\|} \frac{u(x^0+sx)}{s^2}-\frac{\sup_{B_s(x^0)}|u|}{s^2}p \big{\|}_{C^{1,\beta}}-
\big{\|}
\frac{\sup_{B_s(x^0)}|u|}{s^2}p-\frac{\Pi(u(x^0+sx))}{s^2}
\big{\|}_{C^{1,\beta}}
$$
$$
=\frac{\sup_{B_s(x^0)}|u|}{s^2}\bigg(
\big{\|} \frac{u(x^0+sx)}{\sup_{B_s(x^0)}|u|}-p \big{\|}_{C^{1,\beta}}
$$
$$
-\big{\|}
p-\frac{\sup_{B_1(x^0)}|\Pi(u(x^0+sx))|}{\sup_{B_s(x^0)}|u|}\frac{\Pi(u(x^0+sx))}{\sup_{B_1}|\Pi(u(x^0+sx))|}
\big{\|}_{C^{1,\beta}}\bigg).
$$
As a direct consequence of Lemma \ref{normblow} we obtain
$$
\Big| \frac{\sup_{B_1}|\Pi(u(x^0+sx))|}{\sup_{B_s(x^0)}|u|}-1 \Big| \le \frac{C_3 s^2}{\sup_{B_s(x^0)}|u|}.
$$
This, together with Theorem \ref{auto}, implies that
$$
\big{\|}
p-\frac{\sup_{B_1}|\Pi(u(x^0+sx))|}{\sup_{B_s(x^0)}|u|}\frac{\Pi(u(x^0+sx))}{\sup_{B_1}|\Pi(u(x^0+sx))|}
\big{\|}_{C^{1,\beta}}\le C_4 \Big(\frac{s^2}{S^u(x^0,s)}\Big)^\alpha.
$$
Rearranging terms in
(\ref{c1betarotationsetimate}) we get
$$
\Big{\|} \frac{u(x^0+sx)}{\sup_{B_s(x^0)}|u|}-p \Big{\|}_{C^{1,\beta}}\le
C_5\Big( \frac{s^2}{S^u(x^0,s)} \Big)^\alpha\le C(d,M,\alpha,\beta) \big( \frac{\delta}{1+\delta\log(r/s)} \big)^\alpha.
$$
This proves (i).

$ $

Rotating the coordinate system we may assume that
$p^{x^0,u}=p=2x_1x_2$. The first part of the Corollary implies that
$$
u(x^0+s\cdot)<0\textrm{ in } \big\{ (x_1,x_2)\in B_1:\; x_1x_2\le -C(d,M,\alpha,\beta) \Big(\frac{\delta}{1+\delta\log(r/s)} \Big)^{\alpha}\big\}\equiv K_s^-,
$$
that
$$
u(x^0+s\cdot)>0\textrm{ in } \big\{ (x_1,x_2)\in B_1:\; x_1x_2\ge C(d,M,\alpha,\beta) \Big(\frac{\delta}{1+\delta\log(r/s)} \Big)^{\alpha}\big\}\equiv K_s^+
$$
and that
$$\left|\partial_\theta \frac{u(x^0+sx)}{\sup_{B_s(x^0)}|u|}\right|\ge c_6|x|
\textrm{ in } B_1\setminus (K_s^-\cup K_s^+).$$ 
From the implicit function
theorem it follows that, for each $\epsilon>0$, $\{ u=0 \}$ consists of four 
$C^1$-curves in $B_s(x^0)\setminus B_{s/2}(x^0)$. To show that
$\{u=0\}$ consists of two $C^1$-curves we only need to show that these
four curves are differentiable at $x^0$ and that their derivatives match.

The normal $\nu$ of $\{ u=0 \}$ will point in the same (or opposite) direction
as $\nabla u$ at any point of $\big( B_s(x^0)\setminus \{x^0\}\big)\cap \{u=0\}$. Let us consider a point $x^0+sx$ of $\{u=0\}$ such that
$x_2=1$ and  $|x_1|\le 1$: from {\sl (i)} it follows that at the point
$x^0+sx$,
$$
\frac{\nabla \big( u(x^0+sx)\big)}{\sup_{B_s(x^0)}|u|}= \Big(
\frac{\nabla \big(u(x^0+sx)\big)}{\sup_{B_s(x^0)}|u|}-2\nabla (x_1x_2)\Big)+2\nabla(x_1x_2)
$$
$$
=2e_1+ \textrm{ terms of order }
\Big(\frac{\delta}{1+\delta\log(r/s)} \Big)^{\alpha}.
$$
By a similar argument for each of the four components of $\{u=0\}\cap \big( B_s(x^0)\setminus \{x^0\} \big)$
it follows that each component is a $C^1$-curve with modulus of
continuity $\sigma(s)=C_7(\log(r/s))^{-\alpha}$ and that each
component approaches $x^0$ tangentially relative to the $x^1$- or $x^2$-axis. This
proves {\sl (ii)}.\qed

 \bibliographystyle{plain}

\end{document}